\documentclass[12pt,a4paper]{amsart}

\usepackage{amsfonts, amsmath, amssymb, amsthm, amscd, hyperref}

\usepackage{anysize}

\newtheorem{thm}{Theorem}[section]
\newtheorem*{thm*}{Theorem}
\newtheorem{lem}[thm]{Lemma}

\newtheorem{pro}[thm]{Proposition}

\theoremstyle{definition}
\newtheorem{defn}[thm]{Definition}

\theoremstyle{remark}
\newtheorem{rem}[thm]{Remark}

\numberwithin{equation}{section}

\newcommand{\pt}{\mathrm{pt}}

\DeclareMathOperator{\conv}{conv}

\renewcommand{\int}{\mathop{\rm int}}

\renewcommand{\epsilon}{\varepsilon}

\begin{document}

\title[Analogues of the central point theorem\dots]{Analogues of the central point theorem for families with $d$-intersection property in $\mathbb R^d$}

\author{R.N.~Karasev}
\thanks{Supported by the Dynasty Foundation, the President's of Russian Federation grant MD-352.2012.1, the Russian Foundation for Basic Research grants 10-01-00096 and 10-01-00139, the Federal Program ``Scientific and scientific-pedagogical staff of innovative Russia'' 2009--2013, and the Russian government project 11.G34.31.0053.}

\email{r\_n\_karasev@mail.ru}

\address{Roman Karasev, Dept. of Mathematics, Moscow Institute of Physics and Technology, Institutskiy per. 9, Dolgoprudny, Russia 141700}
\address{Roman Karasev, Laboratory of Discrete and Computational Geometry, Yaroslavl' State University, Sovetskaya st. 14, Yaroslavl', Russia 150000}

\keywords{central point theorem, Tverberg's theorem, Helly's theorem}
\subjclass[2000]{52A20, 52A35, 52C35}

\begin{abstract}
In this paper we consider families of compact convex sets in $\mathbb R^d$ such that any subfamily of size at most $d$ has a nonempty intersection. We prove some analogues of the central point theorem and Tverberg's theorem for such families.  
\end{abstract}

\maketitle

\section{Introduction}

Let us start with a definition:

\begin{defn}
A family of sets $\mathcal F$ has \emph{property $\Pi_k$} if for any nonempty $\mathcal G\subseteq \mathcal F$ such that $|\mathcal G|\le k$ the intersection $\bigcap\mathcal G$ is not empty.
\end{defn}

Helly's theorem~\cite{helly1923} states that a finite family of convex sets (or any family of convex compact sets) with $\Pi_{d+1}$ property in $\mathbb R^d$ has a common point. In the review~\cite{eck1993} Helly's theorem and many its generalizations are considered in detail.

In this paper we concentrate on the families with $\Pi_d$ property in $\mathbb R^d$, the ``almost'' Helly property. The typical example of a family with $\Pi_d$ property is any family of affine hyperplanes in general position. It can be easily seen already in the case of affine hyperplanes that such a family need not have a common point, and even need not have a bounded \emph{piercing number}, which is the smallest size of a finite set intersecting any set of the family. The reader may also consult~\cite{kar2008tr}, where some bounds on the piercing number following from the $\Pi_d$ property are given for particular families of sets, for example, balls of equal radii, balls of arbitrary radii, or translates of a single convex compact set in the plane. 

An important consequence of Helly's theorem is the central point theorem~\cite{grun1960,neumann1945,rado1946} for measures: For every absolute continuous probability measure $\mu$ on $\mathbb R^d$ one can find a \emph{central point}, that is a point $x$ such that any halfspace $H\ni x$ has $\mu(H)\ge \frac{1}{d+1}$. Here we discuss the discrete central point theorem for finite point sets instead of measures: For a finite set $X\subset \mathbb R^d$ there exists a \emph{central point} $x\in\mathbb R^d$ such that any half-space $H\ni x$ contains at least 
$$
r = \left\lceil\frac{|X|}{d+1}\right\rceil
$$
points of $X$. Here $|X|$ denotes the cardinality of $X$.

In~\cite{kar2008dcpt} several ``dual'' analogues of the central point theorem were established for the families of affine hyperplanes. For example, if $\mathcal F$ is a family of affine hyperplanes in $\mathbb R^d$ then there exists a point $x\in\mathbb R^d$ such that any ray staring at $x$ intersects at least 
$$
r = \left\lceil \frac{|\mathcal F|}{d+1} \right\rceil
$$
affine hyperplanes of $\mathcal F$. The word ``dual'' here does not mean that this theorem follows from the original discrete central point theorem by the projective duality or some other transformation; this ``dual'' theorem in fact requires a separate proof using some topology.

Here we prove an analogue of the dual central point theorem for every family of convex compact sets with $\Pi_d$ property:

\begin{thm}
\label{pd-cpt} Let a finite family $\mathcal F$ of convex closed sets in $\mathbb R^d$ have property $\Pi_d$. Then there exists a point $x\in \mathbb R^d$ such that any unbounded continuous curve that passes through $x$ intersects at least 
$$
r = \left\lceil \frac{|\mathcal F|}{d+1} \right\rceil
$$
sets in $\mathcal F$.
\end{thm}

Similar to what is done in~\cite{kar2008dcpt} it is natural to generalize this theorem in the spirit of Tverberg's theorem~\cite{tver1966}. First, we have to make a definition. For a family $\mathcal G$ of $d+1$ compact convex sets in $\mathbb R^d$ with $\Pi_d$ property we have two alternatives: either all the sets in $\mathcal G$ have a common point, or the nerve of the family $\mathcal G$ (see~\cite{eck1993} for the discussion of nerves) is a simplicial complex equal to the boundary $\partial \Delta^d$ of the standard $d$-simplex. By the nerve theorem the union $\bigcup \mathcal G$ is homotopy equivalent to $\partial \Delta^d$ (or a $(d-1)$-dimensional sphere) and by the Alexander duality~\cite[Theorem~3.44]{hatcher2002} the complement $\mathbb R^d\setminus\bigcup\mathcal G$ consists of two connected components, one being bounded and the other being unbounded. Now it is natural to make a definition:

\begin{defn}
Consider a family $\mathcal G$ of $d+1$ convex compact sets in $\mathbb R^d$ with $\Pi_d$ property. If the family $\mathcal G$ has no point in common, then the complement of its union consists of two connected components: $X$ and $Y$, where $X$ is bounded and $Y$ is unbounded. In this case for any point $x\in X$ we say that \emph{$\mathcal G$ surrounds $x$}. 
\end{defn}

\begin{rem}
A typical example is: $d+1$ facets of any $d$-dimensional simplex surround any point in the interior of the simplex.
\end{rem}

Now we state the analogue of the Tverberg theorem:

\begin{thm}
\label{pd-dualtver}
Let a finite family $\mathcal F$ of convex compact sets in $\mathbb R^d$ have property $\Pi_d$. Suppose the number
$$
r = \left\lceil \frac{|\mathcal F|}{d+1} \right\rceil
$$
is a prime power. Then there exists a point $x\in \mathbb R^d$ and $r$ pairwise disjoint nonempty subfamilies $\mathcal F_1,\ldots, \mathcal F_r\subseteq\mathcal F$ such that the following condition holds for any $i=1,\ldots, r$:

1) either some member of $\mathcal F_i$ contains $x$;

2) or the family $\mathcal F_i$ surrounds $x$.
\end{thm}

We conjecture that this result holds without the assumption that $r$ is a prime power. In this case this would imply Theorem~\ref{pd-cpt} directly, because any unbounded continuous curve through $x$ must intersect some element of every $\mathcal F_i$. It turns out that in order to deduce Theorem~\ref{pd-cpt}, it is sufficient (see Section~\ref{proof-pd-cpt}) to prove Theorem~\ref{pd-dualtver} only for prime numbers $r$.

It is also possible to give a generalization of Theorem~\ref{pd-dualtver} in the spirit of Tverberg's transversal conjecture~\cite{tvervre1993}; see also~\cite{dol1987,kar2007tt,kar2008dcpt,vre2003,ziv1999} for proofs of some particular cases of Tverberg's transversal conjecture and similar results.

\begin{defn}
Consider a family $\mathcal G$ of $d-m+1$ convex compact sets in $\mathbb R^d$ with $\Pi_{d-m}$ property and an affine $m$-subspace $L$. We say that \emph{$\mathcal G$ surrounds $L$} if $\pi(\mathcal G)$ surrounds the point $\pi(L)$, where $\pi$ is the projection along $L$. 
\end{defn}

\begin{thm}
\label{dualtvertr} 
Suppose that each of $m+1$ families $\mathcal F_i$ $(i=0, \ldots, m)$ of convex compact sets in $\mathbb R^d$ have property $\Pi_{d-m}$. Let the numbers 
$$
r_i = \left\lceil \frac{|\mathcal F_i|}{d-m+1} \right\rceil
$$
be powers of the same prime $p$ and 

a) either $p=2$;

b) or $d-m$ is even;

c) or $m=0$.

Then there exists an affine $m$-subspace $L$ and, for every $i=0,\ldots, m$, some $r_i$ pairwise disjoint nonempty subfamilies $\mathcal F_{i1},\ldots, \mathcal F_{ir_i}\subseteq\mathcal F_i$ such that for any $i=0,\ldots,m$ and $j=1,\ldots, r_i$ the following condition holds:

1) either some member of $\mathcal F_{ij}$ intersects $L$;

2) or the family $\mathcal F_{ij}$ surrounds $L$.
\end{thm}

The case $m=0$ is inserted here to make a unified statement with Theorem~\ref{pd-dualtver}. Actually, in this theorem the sets need not be convex, it is sufficient that all their projections to linear $(d-m)$-subspaces are convex; this property is sometimes called $(d-m)$-convexity.

The assumption that $r_i$ are prime powers is essential in the proof of Theorem~\ref{dualtvertr} since the action of a $p$-torus on the configuration space is required, see Section~\ref{proof-dualtvertr}. Of course, it is natural to conjecture that this restriction is not necessary.

While this paper was considered and reviewed in the journal, another paper~\cite{kar2011} with similar results was published. So the content of this paper has a large intersection with that of~\cite{kar2011}.

{\bf Acknowledgments.}

The author thanks V.L.~Dol'nikov for the discussions that have lead to formulation of these results and the unknown referee for numerous helpful suggestions.

\section{Facts from topology}
\label{eq-cohomology}

We consider topological spaces with continuous (left) action of a finite group $G$ and continuous maps between such spaces that commute with the action of $G$. We call them $G$-spaces and $G$-maps. We mostly consider groups $G=(\mathbb Z_p)^k$ for prime $p$, so-called $p$-tori.

For basic facts about (equivariant) topology and vector bundles the reader is referred to the books~\cite{hsiang1975,mishch1998,milsta1974}. The cohomology is assumed with coefficients in $\mathbb F_p$ ($p$ is the same as in the definition of $G$), we omit the coefficients from notation. Let us start from some standard definitions. In this paper we assume \v{C}ech cohomology, it is safe to make such assumptions in results like Lemma~\ref{del-prod-ind}.

\begin{defn}
Denote by $EG$ the classifying $G$-space, which can be thought of as an infinite join $EG=G*\dots *G*\dots$ with diagonal left $G$-action. Denote $BG=EG/G$. For any $G$-space $X$ denote by $X_G=(X\times EG)/G$, and put (\emph{equivariant cohomology in the sense of Borel}) $H_G^*(X) = H^*(X_G)$. It is easy to verify that for a free $G$-space $X$ the space $X_G$ is homotopy equivalent to $X/G$. 
\end{defn}

Consider the algebra of $G$-equivariant cohomology of the point $A_G = H_G^*(\pt) = H^*(BG)$. For a group $G=(\mathbb Z_p)^k$ the algebra $A_G=H_G^*(\mathbb Z_p)$ has the following structure (see~\cite{hsiang1975}). In the case $p$ odd it has $2k$ multiplicative generators $v_i,u_i$ with dimensions $\dim v_i = 1$ and $\dim u_i = 2$ and relations
$$
v_i^2 = 0,\quad\beta{v_i} = u_i,
$$
where we denote by $\beta(x)$ the Bockstein homomorphism. 

In the case $p=2$ the algebra $A_G$ is the algebra of polynomials of $k$ variables $v_1,\ldots, v_k$ of degree one.

We are going to find the equivariant cohomology of a $G$-space $X$ using the following spectral sequence (see~\cite{hsiang1975,mcc2001}):

\begin{pro}
\label{specseqeq}
The natural fiber bundle $\pi_{X_G} : X_G\to BG$ with fiber $X$ gives the spectral sequence with the $E_2$-term
$$
E_2^{x, y} = H^x(BG; \mathcal H^y(X)),
$$
having a structure of a graded $A_G$-module, and converging to a graded $A_G$-module, associated with the filtration of $H_G^*(X)$.

The system of coefficients $\mathcal H^y(X)$ is obtained from the action of $G = \pi_1(BG)$ on the cohomology $H^y(X)$. The differentials of this spectral sequence are homomorphisms (of corresponding degree) of graded $A_G$-modules.
\end{pro}

This proposition implies the following: If the space $X$ is $(n-1)$-connected then the natural map $A^m_G\to H_G^m(X)$ is injective in dimensions $m \le n$.

Any representation of $G$ can be considered as a vector bundle over the point $\pt$, and it has corresponding characteristic classes in $H_G^*(\pt)$. We need the following lemma, that follows from the results of~\cite[Chapter~III, \S~1]{hsiang1975} (see also~\cite{mm1982,vol1992}):

\begin{lem}
\label{euler-nz}
Let $G=(\mathbb Z_p)^k$, and let $I[G]$ be the subspace of the group algebra $\mathbb R[G]$, consisting of elements
$$
\sum_{g\in G} a_g g,\quad \sum_{g\in G} a_g = 0.
$$
Then the Euler class $e(I[G])\not=0\in A_G$ and is not a divisor of zero in $A_G$. 
\end{lem}

In this lemma the assumption that $G=(\mathbb Z_p)^k$ is essential.

We also need the following folklore fact on the Grassmann variety (see~\cite{dol1987,kar2007tt,ziv1999} for its different applications). Consider the canonical bundle over the Grassmann variety $\gamma : E(\gamma) \to G_d^{d-m}$. In the case $p=2$ we consider the variety of non-oriented linear $(d-m)$-subspaces, and for odd $p$ we consider the variety of oriented subspaces.

\begin{lem}
\label{eulergrass} For the Euler class $e(\gamma)$ modulo $p$ the following holds
$$
e(\gamma)^m\neq 0 \in H^{m(d-m)}(G_d^{d-m}; \mathbb F_p),
$$
if either $p=2$, or $d-m$ is even, or $m=0$. In the latter case we put $e(\gamma)^0 = 1\in H^0(G_d^{d-m}; \mathbb F_p)$ by definition.
\end{lem}

It is hard to locate the place where this lemma was proved for the first time (for example, it follows from Schubert calculus); one particular reference for the proof is~\cite[Lemma~8]{kar2007tt}, where this class in the oriented case is shown to be Poincar\'{e} dual to a set of two points with same signs. In the non-oriented case and mod $2$ cohomology this class is Poincar\'{e} dual to a single point, which is a nontrivial $0$-cycle mod $2$.

\section{Topology of Tverberg's theorem}

In Tverberg's theorem and its topological generalizations (see~\cite{bss1981,vol1996} for example) it is important to consider the configuration space of $r$-tuples of points $x_1,\ldots, x_r\in \Delta^N$ with pairwise disjoint supports. Here $\Delta^N$ is a simplex of dimension $N$. Let us make some definitions, following the book~\cite{mat2003}.

\begin{defn}
Let $K$ be a simplicial complex. Denote by $K_{\Delta}^r$ the subset of the $r$-fold product $K^r$, consisting of the $r$-tuples $(x_1,\ldots, x_r)$ such that every pair $x_i, x_j$ ($i\not=j$) has disjoint supports in $K$. We call $K_{\Delta}^r$ the \emph{$r$-fold deleted product of $K$}.
\end{defn}

\begin{defn}
Let $K$ be a simplicial complex. Denote by $K_{\Delta}^{*r}$ the subset of the $r$-fold join $K^{*r}$, consisting of convex combinations $w_1x_1\oplus\dots\oplus w_rx_r$ such that every pair $x_i, x_j$ ($i\not=j$) with weights $w_i,w_j>0$ has disjoint supports in $K$. We call $K_{\Delta}^{*r}$ the \emph{$r$-fold deleted join of $K$}.
\end{defn}

Note that the deleted join is a simplicial complex again, while the deleted product has no natural simplicial complex structure, although it has some cellular complex structure.

The $r$-fold deleted product of the simplex $\Delta^{(r-1)(d+1)}$ is the natural configuration space in Tverberg's theorem, but sometimes it is simpler to use the deleted join because of the following fact. Denote by $[r]$ the set $\{1, \ldots, r\}$ with the discrete topology, the following lemma is well-known, see~\cite{mat2003} for example.

\begin{lem}
The deleted join of the simplex $(\Delta^N)_\Delta^{*r} = [r]^{*N+1}$ is $(N-1)$-connected.
\end{lem}

If $r$ is a prime power $r=p^k$, then the group $G=(\mathbb Z_p)^k$ can be somehow identified with $[r]$, so a $G$-action on $K_\Delta^r$ and $K_\Delta^{*r}$ by permuting the $r$ factor arises. In this case Proposition~\ref{specseqeq} and the above lemma imply that the natural map $A_G^l\to H_G^l((\Delta^N)_\Delta^{*r})$ is injective in dimensions $l\le N$. We need a similar fact for deleted products, following from the next lemma:

\begin{lem}
\label{del-prod-ind}
Let $r=p^k$, $G=(\mathbb Z_p)^k$, and let $K$ be a simplicial complex. If the natural map $A_G^l\to H_G^l(K_\Delta^{*r})$ is injective for $l\le N$, then the similar map $A_G^l\to H_G^l(K_\Delta^r)$ is injective for $l \le N - r + 1$.  
\end{lem}

\begin{proof}
Define the map $\alpha : K^{*r}\to \mathbb R[G]$ as follows. Let $\alpha$ map a convex combination $w_1x_1\oplus \dots\oplus w_rx_r\in K^{*r}$ to $(w_1, \ldots, w_r)\in\mathbb R^r$, the latter space is identified with $\mathbb R[G]$, if we identify the set $[r]$ with $G$. This map is $G$-equivariant. 

Consider the natural orthogonal projection $\pi : \mathbb R[G]\to I[G]$ (the latter $G$-representation is defined in Lemma~\ref{euler-nz}) and the natural inclusion $\iota : K_\Delta^{*r}\to K^{*r}$. The map $\beta = \pi\circ\alpha\circ\iota : K_\Delta^{*r}\to I[G]$ is $G$-equivariant, and it can be easily checked that 
$$
K_\Delta^r = \{y\in K_\Delta^{*r} : \beta(y) = 0\}.
$$

Now assume the contrary: the image of some nonzero element $\xi\in A_G^l$ is zero in $H_G^l(K_\Delta^r)$ and $l\le N-r+1$. We denote the classes in $A_G$ and their natural images in the equivariant cohomology of $G$-spaces by the same letters if it does not lead to confusion. Put $e(I[G]) = e\in A_G^{r-1}$ for brevity. 

The Euler class of a vector bundle is zero outside the zero set of a section of the bundle. Indeed, if $Z$ is the zero set of a section $s$ of a vector bundle $\nu$ over a space $X$, then the restriction of $\nu$ to $X\setminus Z$ has a nonzero section $s$. Therefore $\nu|_{X\setminus Z}$ has zero Euler class and the needed claim follows from the naturality of the Euler class: the Euler class of the restriction is the restriction of the original Euler class.

So we know that $e$ vanishes in $H_G^{r-1}(K_\Delta^{*r}\setminus K_\Delta^r)$. By the standard property of the cohomology product (often used to estimate the Lusternik--Schnirelman category by the cup-length) we obtain that the class $e\xi$ vanishes over $(K_\Delta^{*r}\setminus K_\Delta^r)\cup K_\Delta^r$, that is over the whole $K_\Delta^{*r}$. By Lemma~\ref{euler-nz} $e\xi\not=0\in A_G^{l+r-1}$, and we come to contradiction with the injectivity condition in the statement of this lemma.
\end{proof}

\section{Proof of Theorem~\ref{pd-dualtver}}

It would be sufficient to prove Theorem~\ref{dualtvertr}, since Theorem~\ref{pd-dualtver} is its particular case. Though we give a separate proof for Theorem~\ref{pd-dualtver} to clarify the exposition. The reasoning in this proof (and the subsequent proofs) is essentially the same as in~\cite{kar2008dcpt,kar2011}.

Consider the simplex $\Delta = \Delta^{n-1}$, along with some identification of its vertices with the members of $\mathcal F$. Take some large enough ball $B\subset\mathbb R^d$, containing all the sets of $\mathcal F$ in its interior. The configuration space that we study is $\Delta^r_\Delta\times B$, denote its elements by $(y_1, y_2, \ldots, y_r, p)$. The points $y_i$ in the simplex $\Delta$ will be considered in barycentric coordinates as functions $y_i : \mathcal F\to \mathbb R^+$ each of them having unit sum of values. The condition that an $r$-tuple $(y_1, \ldots, y_r)$ lies in the deleted product means that the supports of these functions are pairwise disjoin.

Put for brevity $\mathbb R^d = V$. Now let us map our configuration space to $V^r$ by the following rule. Let $\pi_K(p)$ be the metric projection of $p$ to $K\in\mathcal F$ sending every $p$ to the closest to $p$ point in $K$; this map is $1$-Lipschitz and therefore continuous. Put 
$$
f(y_1, y_2, \ldots, y_r, p) = \bigoplus_{i=1}^r \sum_{K\in\mathcal F} y_i(K) (\pi_K(p) - p),
$$ 
This map is evidently continuous and $G$-equivariant, if we identify $V^r$ with $V[G]$ ($V$-valued functions on $G$ with $G$-action by right multiplication by $g^{-1}$).

The map $f$ can be considered as a section of a $G$-equivariant vector bundle $V[G]\times \Delta^r_\Delta\times B \to \Delta^r_\Delta\times B$. This bundle is trivial by definition but the action of $G$ makes it equivariantly nontrivial. The relative Euler class of this section\footnote{
The reader is referred to~\cite{kar2008dcpt} for properties of the relative Euler class. It is important that the relative Euler class depends on both the vector bundle and its section.
}
can be decomposed according to the decomposition $V[G] = V\oplus V\otimes I[G]$, multiplicativity of the relative Euler class (see~\cite{kar2008dcpt}), and the K\"unneth formula:
$$
e(f) = w^d\times u\in H_G^{rd}(\Delta^r_\Delta\times B, \Delta^r_\Delta\times \partial B)=H^{d(r-1)}_G(\Delta^r_\Delta) \otimes H^d(B, \partial B).
$$
Here $w$ is the image of $e(I[G])$ in $H^{r-1}(\Delta^r_\Delta)$ and $u$ is the generator of $H^d(B, \partial B)$. We indeed obtain $u$ as the second factor because for any fixed $(y_1,\ldots, y_r)$ the corresponding projection $f''$ of the section $f$ to the summand $V$ corresponds to a vector pointing from $p$ to a convex combination of vectors $\pi_K(p)$. If $B$ contains all sets of $\mathbb F$ in its interior as we have assumed, then this vector always points inside $B$ for $p\in\partial B$. Hence the corresponding Euler class is the same as in the Brouwer fixed point theorem~\cite{br1910}, which is the generator of $H^d(B, \partial B)$.

By Lemmas~\ref{euler-nz} and \ref{del-prod-ind}, $w^d\not=0\in H_G^{d(r-1)}(\Delta^r_\Delta)$, and the K\"unneth formula implies that $e(f)\not=0$. 

The map $f$ therefore must have a zero, let it be $(y_1, y_2, \ldots, y_r, p)$. For any $K\in\mathcal F$ there is at most one $i\in[r]$ such that $y_i(K) > 0$, since $y_i$'s have disjoint supports. In this case we put $K$ to the subset $\mathcal F_i$. From the definition of $f$ it follows that for any $i$ the projections of $p$ to the sets $K\in\mathcal F_i$ have $p$ in their convex hull.

We use the following lemma:

\begin{lem}
\label{surround}
Let a family $\mathcal G = \{G_1, \ldots, G_{d+1}\}$ of convex compact sets in $\mathbb R^d$ have property $\Pi_d$. Let a point $p\in \mathbb R^d$ be such that $p$ lies in the interior of the convex hull of $g_1, \ldots, g_{d+1}$, where $g_i$ is the closest to $p$ point in $G_i$. Then $\mathcal G$ surrounds $p$.
\end{lem}

\begin{proof}[Proof of Lemma~\ref{surround}]
Consider the half-spaces
$$
H_i = \{x\in \mathbb R^d : (x, g_i - p) \ge (g_i, g_i - p)\}
$$
and note that $G_i\subseteq H_i$. Clearly, $\bigcap_{i=1}^{d+1} H_i = \emptyset$.

For any $i=1,\ldots, d+1$ the nonempty intersection $\bigcap_{j\not=i} G_j$ is contained in $\bigcap_{j\not=i} H_i$, take one point $x_i\in \bigcap_{j\not=i} G_j$. The simplex $\Delta=\conv_{i=1}^{d+1} \{x_i\}$ contains $\mathbb R^d\setminus\bigcup_{i=1}^{d+1} H_i\ni p$ (compare~\cite[Lemma~1]{kar2008tr}), and every its facet $\partial_i\Delta=\conv_{j\not=i} \{x_i\}$ is contained in the corresponding $G_i$. 

Thus $p\not\in \bigcup_{i=1}^{d+1} G_i$ and is separated from infinity by $\bigcup_{i=1}^{d+1} G_i\supseteq \partial \Delta$, so $\mathcal G$ surrounds $p$ by definition.
\end{proof}

If $p$ coincides with one of $\pi_K(p)$ (for $K\in\mathcal F_i$), then $p$ is already contained in $\bigcup\mathcal F_i$. If $p$ lies in the interior of the convex hull of some $d+1$ points of $\{\pi_K(p)\}_{K\in\mathcal F_i}$, we reduce $\mathcal F_i$ so that it contains only those $d+1$ corresponding sets $K$ and note that $\{\pi_K(p)\}_{K\in\mathcal F_i}$ surround $p$ by Lemma~\ref{surround}, and therefore $\mathcal F_i$ surrounds $p$.

If none of the above two alternatives holds, then $p$ lies in the relative interior of the convex hull of some $n<d+1$ points $\pi_{K_1}(p),\ldots, \pi_{K_n}(p)$, $K_1,\ldots, K_n\in\mathcal F_i$. Denote the half-spaces 
$$
H_K = \{x\in \mathbb R^d : (x, \pi_K(p) - p) \ge (\pi_K(p), \pi_K(p) - p)\}.
$$
Note that $K\subseteq H_K$ (since $\pi$ is the projection) and the half-spaces $H_{K_1}, \ldots, H_{K_n}$ have empty intersection. So some $n < d+1$ sets of $\mathcal F_i$ have an empty intersection that contradicts the $\Pi_d$ property. So the case $n<d+1$ is impossible and the proof is complete.

\section{Proof of Theorem~\ref{dualtvertr}}
\label{proof-dualtvertr}

For any affine $m$-subspace $L$ denote the unique linear $(d-m)$-subspace in $\mathbb R^d$, orthogonal to $L$, by $L^\perp$. It is easy to see that $L$ is determined uniquely by $L^\perp$ and the point $L\cap L^\perp$. So the variety of all affine $m$-subspaces is the total space of the canonical bundle $\gamma_d^{d-m}$ over the Grassmann variety $G_d^{d-m}$.

For any $V \in G_d^{d-m}$ denote the orthogonal projection of $\mathbb R^d$ onto $V$ by $\pi_V$. For any $X\in \bigcup_{i=0}^m\mathcal F_i$, $V\in G_d^{d-m}$, and $p\in V$ denote by $\phi(V, p, X)$ the closest to $p$ point in $\pi_V(K)$. This point depends continuously on the pair $(V, p)$ (a standard technical argument showing this is omitted) and lies in $V$.

Fix some $i=0,\ldots, m$ and define a linear map $\psi_i: K_i = \Delta^{|\mathcal F_i|+1} \to V$ that maps the vertices of the simplex (corresponding to $X\in \mathcal F_i$) to the points $\phi(V, p, X)-p$ for $X\in\mathcal F_i$ and is piece-wise linear. Denote by $\xi_i : (K_i)_\Delta^{r_i}\to V^{r_i}$ the corresponding map of the deleted product. This map is the analogue of $f$ from the previous proof, but we have to define one such map for every $\mathcal F_i$.

Let the group $G_i=(\mathbb Z_p)^{k_i}$, where $r_i=p^{k_i}$, act on the deleted product $L_i = (K_i)_\Delta^{r_i}$ and on $V^{r_i}$ by permutations, we put $V^{r_i} = V[G_i]$ to indicate this action, the map $\xi_i$ thus becomes $G_i$-equivariant.

In the sequel we put $\gamma_d^{d-m}=\gamma$ for brevity. Summing up all the maps we obtain a map
$$
\xi : L_0\times\dots\times L_m\to V[G_0]\oplus\dots\oplus V[G_m].
$$
The map $\xi$ also depends on the pair $(V, p)\in E(\gamma)$ continuously, so actually it gives a section $\xi$ of the vector bundle $\nu$ with fiber $V[G_0]\oplus\dots\oplus V[G_m]$ over the space $E(\gamma) \times L_0\times\dots\times L_m$. Here $V$ as a function of the pair $(V, p)$ can be treated as the pullback of the vector bundle $\gamma\to G_d^{d-m}$ by the map $\gamma : E(\gamma)\to G_d^{d-m}$. We denote this pullback by $\gamma$ (it does not make a confusion) and therefore assume $\xi$ to be a section for $\gamma \otimes (\mathbb R[G_0]\oplus\dots\oplus \mathbb R[G_m])$.

To prove the theorem we have to find $V\in G_d^{d-m}, p\in V, (y_0, \ldots, y_m)\in L_0\times\dots\times L_m$ such that $\xi(V, p, y_0, \ldots, y_m) = 0$.

If we take the bundle of large enough balls $B(\gamma)$ in $\gamma$, the section $\xi$ obviously has no zeros over $\partial B(\gamma)\times L_0\times\dots\times L_m$ (this happens when all the balls $\partial B(\gamma)$ contain all the projections $\pi_V(X)$ for $X\in\bigcup_i \mathcal F_i$ in their interiors). To guarantee the zeros for the section $\xi$, we have to find the relative Euler class
$$
e(\xi)\in H_{G_0\times\dots\times G_m}^{(d-m)(r_0+\dots+r_m)} (B(\gamma)\times L_0\times\dots\times L_m, \partial B(\gamma)\times L_0\times\dots\times L_m).
$$ 
Put for brevity $G=G_0\times\dots\times G_m$.

Let us decompose the bundle $\nu$ and its section $\xi$ in the following way. Any $V[G_i]$ can be split $V[G_i] = V\otimes \mathbb R[G] = V\otimes\mathbb R\oplus V\otimes I[G_i]=V\oplus V\otimes I[G_i]$. So the section $\xi$ splits into section $\eta$ of the bundle $\upsilon=\gamma^{m+1}$ and $\zeta$ of the bundle $\omega=\gamma\otimes \bigoplus_{i=0}^m I[G_i]$, and $\nu=\upsilon\oplus\omega$.  

The section $\eta$ has no zeroes on $\partial B(\gamma)\times L_0\times\dots\times L_m$ and, in fact, for large enough balls in $B(\gamma)$ the homotopy $\eta_t = (1-t)\eta + t (-p,\ldots, -p)$ connects it to the section $(-p, \ldots, -p)$ so that $\eta_t$ has no zeroes on $\partial B(\gamma)\times L_0\times\dots\times L_m$ for all $t\in[0,1]$. The section $\eta_1$ does not depend on the factor $L_0\times\dots\times L_m$ and it can be easily seen that (see~\cite{kar2008dcpt}, the proof of Theorem~6)
\begin{multline*}
e(\eta) = u(\gamma)e(\gamma)^m\times 1\in H^{(d-m)(m+1)}(B(\gamma), \partial B(\gamma))\times H_G^0(L_0\times\dots\times L_m)\subset\\
\subset H^{(d-m)(m+1)}(B(\gamma)\times L_0\times\dots\times L_m, \partial B(\gamma)\times L_0\times\dots\times L_m),
\end{multline*}   
where $u(\gamma)$ is the Thom's class of $\gamma$ (in $H^*(E(\gamma))$), $e(\gamma)$ is its Euler class, and the last inclusion is the K\"unneth formula. Lemma~\ref{euler-nz} and the Thom isomorphism show that $u(\gamma)e(\gamma)^m\not = 0$ (compare~\cite[Proof of Theorem~6]{kar2008dcpt}).

Now we consider the class $e(\zeta)\in H_G^{(d-m)(r_0+\dots+r_m-m-1)}(B(\gamma)\times L_0\times\dots\times L_m)$. Taking some fixed $p\in B(\gamma)$ and considering the inclusion 
$$
\iota_p : L_0\times\dots\times L_m=\{p\}\times L_0\times\dots\times L_m\to B(\gamma)\times L_0\times\dots\times L_m
$$
and the induced bundle $\iota_p^* (\omega) = \bigoplus_{i=0}^m (I[G_i])^{d-m}$, we obtain
\begin{multline*}
\iota_p^*(e(\zeta)) = e(I[G_0])^{d-m}\times e(I[G_1])^{d-m}\times \dots \times e(I[G_m])^{d-m}\in H_G^*(L_0\times\dots\times L_m) =\\
= H_{G_0}^*(L_0)\times\dots\times H_{G_m}^*(L_m),
\end{multline*}
the last equality being the K\"unneth formula. By Lemmas~\ref{euler-nz} and \ref{del-prod-ind}, for any $i=0,\ldots, m$ the Euler class $e(I[G])^{d-m}\not=0\in H_{G_i}^{(d-m)(r_i-1)}(L_i)$ and, by the K\"unneth formula, $\iota_p^*(e(\zeta)) = a\not=0$. From one more K\"unneth formula for the product $B(\gamma)\times L_0\times \dots\times L_m$ it follows that 
$$
e(\zeta) = 1\times a + \sum_j b_j\times c_j,
$$ 
where $b_j\in H^*(B(\gamma)), c_j\in H_G^*(L_0\times \dots\times L_m)$ are some classes such that $\dim b_j > 0$ for all $j$. Hence
$$
e(\xi) = u(\gamma)e(\gamma)^m\times a + \sum_j u(\gamma)e(\gamma)^m b_j\times c_j,
$$
and $e(\xi)\not=0$ by the K\"unneth formula (its first summand cannot be eliminated by the latter sum).
 
Now we have a zero of $\xi$ at $(V, p, y_0,\ldots, y_m)$. Every point $y_i\in L_i$ is actually an $r_i$-tuple of points $y_{i1}, \ldots, y_{ir_i}\in K_i = \Delta^{|\mathcal F_i|+1}$ with pairwise disjoint supports. We identify the vertices of $K_i$ with $\mathcal F_i$ and write
$$
y_{ij} = \sum_{X\in \mathcal F_i} w(i, j, X) X.
$$  
Denote $\mathcal F_{ij} = \{X\in\mathcal F_i : w(i, j, X) > 0\}$, each $X$ is assigned to no more than one of $\mathcal F_{ij}$, because $y_{ij}$ have pairwise disjoint supports. The condition $\xi=0$ implies that for any $i=0,\ldots, m$ and $j = 1,\ldots, r_i$ the point $p$ is a convex combination of its projections to the sets $\pi_V(X)$:
$$
p = \sum_{X\in \mathcal F_{ij}} w(i, j, X) \phi(V, p, X).
$$

Now we define $L$ to be the affine subspace, orthogonal to $V$ and passing through $p$. The rest of the proof is the same as in the previous section, because every $\mathcal F_{ij}$ either intersects $L$ (equivalently, the family $\{\pi_V(X)\}_{X\in\mathcal F_{ij}}$ covers $p$) or surrounds $L$ (equivalently, the family $\{\pi_V(X)\}_{X\in\mathcal F_{ij}}$ surrounds $p$).

\section{Proof of Theorem~\ref{pd-cpt}}
\label{proof-pd-cpt}

In this theorem we can assume that $\mathcal F$ consists of compact sets. Indeed, for a large enough ball $B$ the family $\{X\cap B\}_{X\in\mathcal F}$ consists of compact sets and has property $\Pi_d$.

As it was already noted, this theorem follows from Theorem~\ref{pd-dualtver} directly when $r$ is a prime power. Consider some other $r$. Obviously, it is sufficient to prove the theorem in the case $N = |\mathcal F| = (d+1)(r-1)+1$. 

By the Dirichlet theorem on arithmetic progressions, we can find a positive integer $k$ such that $R=k(r-1) + 1$ is a prime. Now take the family $\mathcal F'$ of size $kN$ by simply repeating each set in $\mathcal F$ exactly $k$ times. Note that
$$
kN = k(d+1)(r-1) + k = (d+1)(R-1) + k \ge (d+1)(R-1) + 1,
$$
so we can apply the case of the theorem, that is already proved, to $\mathcal F'$ to get some point $x$.

Every unbounded closed curve $C\ni x$ intersects at least $R=k(r-1)+1$ sets of $\mathcal F'$. Each set of $\mathcal F$ is counted no more that $k$ times, then we conclude that $C$ intersects at least $r$ sets of $\mathcal F$.

\end{document}